\documentclass[a4paper,12pt]{amsart}
\usepackage{amsfonts,amsmath,amssymb,amscd,latexsym}

\newcommand{\bi}{\bibitem}
\newcommand{\nb}{\newblock}

\newcommand{\be}[1]{\begin{equation}\label{#1}}
\newcommand{\ee}{\end{equation}}
\newcommand{\la}{\langle\,}
\newcommand{\ra}{\,\rangle}
\newcommand{\ccc}{{\mathcal C}}
\newcommand{\zz}{{\mathbb Z}}
\newcommand{\fff}{{\mathbb F}}
\newcommand{\rrr}{{\mathbb R}}
\newcommand{\infl}{\mathop{\mathrm{infl}}}
\newcommand{\Sol}{\mathop{\mathrm{Sol}}}
\newcommand{\supp}{\mathop{\mathrm{supp}}}
\newcommand{\Card}{\mathop{\mathrm{Card}}}
\newcommand{\ovr}{\overline}

\newtheorem{thm}{\quad Theorem}
\newtheorem{lm}{\quad Lemma}

\begin{document}

\title{Strict dead end elements in free soluble groups}
\author{Victor Guba}
\address
{Vologda State Pedagogical University,
6 S.\,Orlov Street, Vologda 160600 Russia}
\email
{guba{@}uni-vologda.ac.ru}
\thanks{This research is partially supported by the RFFI
grant 05--01--00895 and the INTAS grant 99--1224.}
\subjclass[2000]{Primary 20F32; Secondary 05C25}
\keywords{Cayley graphs; dead end elements; free soluble groups}

\begin{abstract}

Let $G$ be a group generated by a finite set $A$. An element $g\in G$
is a strict dead end of depth $k$ (with respect to $A$) if
$|g|>|ga_1|>|ga_1a_2|>\cdots>|ga_1a_2\cdots a_k|$ for any $a_1,a_2,
\ldots, a_k\in A^{\pm1}$ such that the word $a_1a_2\cdots a_k$ is freely
irreducible. (Here $|g|$ is the distance from $g$ to the identity in
the Cayley graph of $G$.) We show that in finitely generated free soluble
groups of degree $d\ge2$ there exist strict dead elements of depth
$k=k(d)$, which grows exponentially with respect to $d$.

\end{abstract}

\maketitle

Let $G$ be a group generated by a finite set $A$. By $\ccc=\ccc(G,A)$
we denote the {\em right Cayley graph\/} of $G$ with respect to $A$.
The set of vertices of this graph is $G$, the set of edges is
$G\times A^{\pm1}$. Every edge $e=(g,a)$, where $g\in G$, $a\in A^{\pm1}$,
starts in $g$ and ends in $ga$. This edge is labelled by $a$. An inverse
edge is $e^{-1}=(ga,a^{-1})$. For every $g\in G$ let $|g|$ denote the
distance in $\ccc$ from the vertex $g$ to the identity.

An element $g$ is a {\em dead end\/} in $\ccc$ whenever $|g|>|ga|$
for all $a\in A^{\pm1}$. This concept was introduced by Bogopolski
in \cite{Bog}. One of the easiest examples is the following. Let
us consider the group $G=\zz\times\zz_2$, where $a$, $b$ generate the
factors. It is easy to see that the Cayley graph of $G$ in generators
$a,b$ has no dead end vertices. However, if the generating set is
$A=\{a,ab\}$, then $b$ becomes a dead end element. This example
shows that being a dead end element strongly depends on the
generating set. However, if the generating set is clear, we will
often say that an element of $G$ is (or is not) a dead end.

A less trivial example of dead end elements comes from the
R.\,Thomp\-son group $F$. Recall that $F$ is the group presented by
$$
\la x_0,x_1,x_2,\ldots\mid x_jx_i=x_ix_{j+1}\ (j>i)\ra.
$$
This group can be generated by $x_0$, $x_1$. In these generators,
the group has infinitely many dead end elements. Say,
$g=x_0^2x_1x_6x_3^{-1}x_0^{-2}$ will be a dead end element of length $11$.
One can see details in \cite{Gu04}. (There are no dead end elements in $F$
of length less than $11$.) There are descriptions of all dead end elements
in $F$ given independently in \cite{BeBr} and \cite{ClTa}. They are based
on the representations of $F$ in terms of (pairs of) rooted binary forests
and trees, respectively.

One can try to generalize the concept of a dead end element in different
ways. Say, one can ask if there are elements in $F$ with the following
property: $|g|>|ga|>|gab|$ for all $a,b\in\{x_0^{\pm1},x_1^{\pm1}\}$
such that $a$, $b$ are not inverses of each other. It can be shown that
such elements do not exist \cite{ClTa}.

One of the generalizations has been considered in \cite{ClTa2}, see
also \cite{ClR}. Let $G$ be a group and let $A$ be its finite generating
set. An element $g\in G$ is a {\em dead end of depth\/} $k$ in the Cayley
graph $\ccc(G,A)$ whenever the ball of radius $k$ around $g$ is
contained in the ball of radius $|g|$ around the identity. In other
words, $|g|\ge|ga_1|,|ga_1a_2|,|ga_1a_2\cdots a_k|$ for all
$a_1,a_2,\ldots, a_k\in A^{\pm1}$. Hence each dead end element is
a dead end of depth $2$. It follows from \cite{ClTa} that $F$ has
no dead ends of depth $\ge3$ (the set of generators is $\{x_0,x_1\}$
as above).

One can ask whether there are groups whose Cayley graphs contain
dead end elements of arbitrary depth. A positive answer is given
in \cite{ClTa2}, where the group is the {\em lamplighter group\/},
which is the restricted wreath product $\zz_2\wr\zz$ with natural
set of generators. Finitely presented examples of that kind are
considered in \cite{ClR}.
\vspace{1ex}

One more natural generalization comes as follows. Let us say that
an element $g\in G$ is a {\em strict dead end of depth\/} $k$ in the
Cayley graph $\ccc(G,A)$ whenever
$|g|>|ga_1|>|ga_1a_2|>\cdots>|ga_1a_2\cdots a_k|$ for all
$a_1,a_2,\ldots a_k\in A^{\pm1}$ such that the word $a_1a_2\cdots a_k$
is freely irreducible. This means that if we travel in $\ccc$ from $g$
along a path of length $k$ with freely irreducible label, then we
become closer and closer to the identity each time we go through an
edge. Obviously, each dead end element is exactly a strict dead end
element of depth $1$.

Unlike the case of (non-strict) dead end elements, it is obvious
that any group can have strict dead elements of bounded depth only.
Indeed, a group freely generated by $A$ has no dead ends at all. If
there is a nontrivial group relation between elements of $A$, then
let us take the shortest one. If $\rho$ denotes its length, then there
are no strict dead end elements of depth $k>[\,\rho/2\,]$. The goal
of this article is, given a $k$, to construct a group that has strict
dead elements of depth $k$. Moreover, we find such examples among
free soluble groups. These groups are not finitely presented if the
degree of solvability is greater than or equal to $2$ (see \cite{Shm}).
However, an easy observation shows that the needed examples can be
found in the class of finitely presented groups and even finite groups.
\vspace{1ex}

We need to recall some facts from the Stallings theory of group ends.
Let $G$ be a finitely generated group. Let us consider its Cayley
graph $\ccc=\ccc(G,A)$ and take the ball $B_n$ of radius $n$ around the
identity. The graph $\ccc-B_n$ (we delete all vertices at distance
$\le n$ from the identity and all edges coming out of them) may consist
of several connected components. If $n$ approaches infinity, then the
number of infinite connected components of $\ccc-B_n$ approaches some
value that should be equal to $0$, $1$, $2$, or infinity. This is the
number of {\em ends\/} of the group $G$. It is independent on the choice
of the finite generating set for $G$. It is obvious that $G$ has zero ends
if and only if it is finite. If $G$ has two ends, then it is a finite
extension of $\zz$. If $G$ has infinitely many ends, then $G$ is a
nontrivial amalgamated product or a nontrivial HNN-extension with finite
amalgamated (associated) subgroup(s) \cite{Sta}. Therefore, in many cases
we can conclude that a group has only one end.
\vspace{1ex}

Now we need to recall some facts from \cite{DLS} that will play an
important r\^ole for us. For some of them, we recall the idea of the
proof for reader's convenience.

Let $\fff_m$ be a free group with $m$ generators and let $R$ be its
normal subgroup. We fix the standard generating set $A$ for $\fff_m$.
This will also denote the generating sets for quotient groups $\fff_m/R$
and $\fff_m/R'$, where $R'=[R,R]$ is the commutator subgroup of $R$.

By $\Sol(m,d)$ we denote the {\em free soluble\/} group of rank $m$
and {\em degree of solvability\/} $d$. This is the quotient
$\fff_m/\fff_m^{(d)}$, where $H^{(d)}$ denotes the $d$th {\em derived
subgroup\/} of $H$. This is defined inductively: $H^{(0)}=H$ and
$H^{(n+1)}=[H^{(n)},H^{(n)}]$ for all $n\ge0$. To work with free
soluble groups, one needs to describe the connection between groups
$\fff_m/R$ and $\fff_m/R'$ in its general form.

There is a natural way to describe the solution to the word problem
in $\fff_m/R'$ based on some properties of the Cayley graph of
$\fff_m/R$. Moreover, the paper \cite{DLS} contains an algorithm
how to find the length of a given element in $\fff_m/R'$ with respect
to $A$. We need to introduce some terminology.
\vspace{1ex}

Let $\Gamma$ be a graph in the sense of Serre \cite{Serre} (that is,
each edge has an inverse). A {\em flow\/} on $\Gamma$ is a function
$\mu$ from the set of edges $E=E(\Gamma)$ to $\rrr$ such that
$\mu(e^{-1})=-\mu(e)$ for all $e\in E$. We say that $\mu(e)$ is a flow
{\em through\/} the edge $e$. If $\Gamma$ is locally finite, then we
can define an {\em inflow\/} $\infl_\mu(v)$ into a vertex $v$. This is
the sum of flows through all edges $e$ that have $v$ as its terminal
vertex.

Given a flow $\mu$ on a graph, we will say that its {\em support\/}
$\supp(\mu)$ is a subgraph in $\Gamma$ that is the union of all the
edges $e$ satisfying $\mu(e)\ne0$ together with their endpoints.
Thus $\mu$ has finite support if and only if $\mu(e)=0$ for all but
finitely many $e\in E$. For an edge $e\in E$, we denote by $\chi_e$
the {\em characteristic flow\/} of $e$ on $\Gamma$. Namely, $\chi_e(e)=1$,
$\chi_e(e^{-1})=-1$, $\chi_e(x)=0$ otherwise. Obviously, the set of flows
on $\Gamma$ forms a vector space. If we fix an orientation on $\Gamma$
(that is, if we choose a subset $E_+$ in $E$ such that $E_+$ contains
exactly one of the edges $e$, $e^{-1}$ for all $e\in E$), then the set of
flows $\chi_e$, $e\in E_+$ forms a basis of the subspace of flows with
finite support.

A flow on $\Gamma$ is called {\em balanced\/} whenever the inflow into
any vertex is zero. A flow $\mu$ on $\Gamma$ is {\em semi-balanced\/}
whenever there exist two vertices $v_-$, $v_+$ such that
$\infl_\mu(v_-)=-1$, $\infl_\mu(v_+)=1$ and $\infl_\mu(x)=0$ for all
$x\ne v_-,v_+$. Every path $p$ in $\Gamma$ induces a flow $\chi_p$ on
$\Gamma$ in the following way: if $p=e_1\cdots e_n$ is a product of edges,
then $\chi_p=\chi_{e_1}+\cdots+\chi_{e_n}$ by definition. Such a flow is
always balanced or semi-balanced.

Given a flow $\mu$, an edge is called a $P$-{\em edge\/} whenever
$\mu(e)>0$. If $\supp(\mu)$ is finite, then the sum of flows through
all the $P$-edges is called the {\em weight\/} of $\mu$.

\begin{lm}
\label{connflow}
Let $\mu$ a nonzero integer-valued flow on a graph $\Gamma$. Suppose
that $\mu$ is balanced and let $\supp(\mu)$ be finite and connected.
Then for any vertex $v$ of $\supp(\mu)$, there exists a loop $p$ at
$v$ in $\Gamma$ such that $p$ is a product of $P$-edges and every edge
$e\in E$ occurs in this product exactly $\mu(e)$ times.
\end{lm}

\proof
Let $N$ be the weight of $\mu$. We proceed by induction on $N$. Let
$e_1$ be any $P$-edge. For every $i\ge1$, let $v_i$ denote the terminal
vertex of $e_i$. Since the inflow into $v_i$ is zero, there exists a
$P$-edge $e_{i+1}$ starting at $v_i$. Since the support of $\mu$ is finite,
some edges of the form $e_i$ ($i\ge1$) will repeat. Hence we can find a
closed path of the form $q=e_{i+1}\cdots e_{i+k}$, $k\ge1$ that consists
of different $P$-edges.

Let $\mu'=\mu-\chi_q$. Clearly, $\mu'$ is also balanced and has a finite
support. The weight of $\mu'$ equals $N-k<N$. If $\mu'$ is zero, then we
are done. Otherwise $\supp(\mu')$ is a union of its connected components
$\Sigma_1$, \dots, $\Sigma_s$ ($s\ge1$). We can apply the inductive
assumption to the restriction of $\mu'$ on each of the components.

It is easy to see that each $\Sigma_j$ ($1\le j\le s$) has a common
vertex $v_j$ with the path $q$. These vertices are different for
different values of $j$. Changing the order of the components, we can
assume that $q=q_0q_1\cdots q_s$, where $v_j$ is the initial point of
$q_j$ ($1\le j\le s$). By the inductive assumption, there is a loop $p_j$
at $v_j$ that contains each $P$-edge $e$ of $\Sigma_j$ exactly $\mu'(e)$
times ($1\le j\le s$). Now the path $p=q_0p_1q_1\cdots p_sq_s$ is the
loop we wanted to find.
\endproof

Let $G$ be a group generated by a finite set $A$. For every group word
$w$ over $A$ there exists a unique path $p=p(w)$ in the Cayley graph
$\ccc=\ccc(G,A)$ starting at the identity and labelled by $w$. So we
can assign to $w$ the flow $\mu_{p(w)}$ on $\ccc$. This flow is called
the {\em flow induced by\/} the word $w$. It is obvious that if two
words are equal in the free group, then they induce the same flow. It
is also obvious that the flow induced by a word has finite support. It
is balanced whenever $p=p(w)$ is closed (that is, $w$ equals $1$ in $G$)
and semi-balanced whenever $p=p(w)$ is not closed (in this case $v_-=1$,
$v_+=g$, where $g\in G$ is the element represented by $w$). It is easy
to prove the following fact.

\begin{lm}
\label{bal}
A flow $\mu$ on $\ccc=\ccc(G,A)$ is induced by a word if and only if
the following conditions hold:
\begin{itemize}
\item
$\mu$ is integer-valued,
\item
$\mu$ has finite support,
\item
$\mu$ is balanced or it is semi-balanced with $v_-=1$.
\end{itemize}
\end{lm}

\proof
The ``only if" part is obvious. Now assume that a flow $\mu$ satisfies
all the three conditions. First of all, let $\mu$ be balanced. If
$\supp(\mu)$ is connected, we find a path $q$ from the identity to
a point $v$ in $\supp(\mu)$ and a loop $p$ from Lemma \ref{connflow}.
Then the label of $qpq^{-1}$ will induce the desired flow. If $\supp(\mu)$
is not connected, then we take a product of words that induce the flow
for each of the connected components.

In the semi-balanced case, let us assume that $\mu$ has the property
$\infl_\mu(1)=-1$, $\infl_\mu(g)=1$ with zero inflow into all the other
vertices. Let us consider any path $q$ labelled by a word $Q$ that
connects $g$ and $1$. Let us add the flow $\chi_q$ to $\mu$. This flow
will be balanced and thus induced by a word $w$. Clearly, the word
$wQ^{-1}$ will induce $\mu$.
\endproof

The following fact has been proved in \cite{DLS}.

\begin{lm}
\label{wp}
Let $R$ be a normal subgroup in the free group $\fff_m$. Two words
$w_1$, $w_2$ are equal modulo the commutator subgroup $R'$ if and
only if they induce the same flow on the Cayley graph of $\fff_m/R$
with respect to the natural set of generators.
\end{lm}

For the sake of completeness, let us recall the idea of the proof.
A word $w$ belongs to $R$ if and only if the path $p=p(w)$ is a loop
in the Cayley graph $\ccc$ of $\fff_m/R$. Hence any word of the form
$[r_1,r_2]=r_1^{-1}r_2^{-1}r_1r_2$, where $r_1,r_2\in R$, will induce
a zero flow on $\ccc$. The same holds for a product of these commutators,
which implies the ``only if" part. To prove the converse, it suffices to
choose a spanning subtree $T$ in $\ccc$. For each edge $e$, one can
define a loop $\ovr{e}=peq^{-1}$ at the identity, where $p$, $q$ are
geodesic paths in $T$ that connect $1$ with the initial and the terminal
vertex of $e$, respectively. If $w_1$, $w_2$ induce the same flow on
$\ccc$, then the paths $p(w_1)$, $p(w_2)$ have the same endpoints. Hence
$s=p(w_1)p(w_2)^{-1}$ will be the loop at the identity. If
$s=e_1\cdots e_n$ is the product of edges, then $s$ is freely homotopic
to the path $\bar e_1\cdots\bar e_n$. Since the flow $\chi_s$ is zero,
the number of occurrences of $e$ in $s$ is the same as the number of
occurrences of $e^{-1}$, for any $e$. Taking into account that $\bar e$
represents an element of $R$ and the fact that $\ovr{e^{-1}}=\ovr{e}^{-1}$,
we conclude that the label of $s$ equals $1$ modulo $R'$ (we do
permutations in the product $\ovr{e_1}\cdots\ovr{e_n}$ cancelling each
$\ovr{e}$ with its inverse).
\vspace{2ex}

The paper \cite{DLS} also has a description of lengths of all elements in
$\fff_m/R'$ and their minimal word representatives. This means that, given
a word $w$, we can find a word $v$ that equals $w$ modulo $R'$ and has the
shortest possible length. (This word is not unique in general.) The length
of the shortest representative of $g$ is denoted by $|g|$ and it is called
the {\em length\/} of the element $g\in\fff_m/R'$. In fact, the algorithm
to find the length of an element given by a word $w$ is polynomial on
$|w|$. (One can compare this situation with the one described in
\cite{Par}, where it is shown that the problem to find the length of an
element of the restricted wreath product $\zz\wr\zz^2$ in its natural
generators turns out to be NP-complete.)

We need an explicit description of lengths of elements in $\fff_m/R'$
of some particular form.

Lemma \ref{bal} shows us that any balanced flow with finite support on
the Cayley graph of $\fff_m/R$ is induced by a word $w$. It follows
from Lemma \ref{wp} that $w$ is defined uniquely modulo $R'$. Therefore,
we can say that any flow with the above properties defines canonically
an element of $\fff_m/R'$. (In fact, all these elements form the
subgroup $R/R'$.)

The following statement is a partial case of \cite[Theorem 2]{DLS}.

\begin{lm}
\label{lng}
Let $\mu$ be a nonzero balanced flow with finite support on the Cayley
graph of $\fff_m/R$. Let $g\in\fff_m/R'$ be defined by this flow. Suppose
that $\supp(\mu)$ is connected. Let $d$ be the distance from the identity
to $\supp(\mu)$ in the above Cayley graph. Then the length of the shortest
word that represents $g$ in $\fff_m/R'$ is given by the formula $|g|=N+2d$,
where $N$ is the weight of $\mu$.
\end{lm}

\proof
Suppose that a word $w$ induces the flow from the statement. Let $p=p(w)$
be the path labelled by $w$ starting at the identity in the Cayley graph
$\ccc$ of the group $\fff_m/R$. Let $q$ be the shortest subpath in $p$
that contains all edges with nonzero flow. Obviously, $|q|\ge N$. If
$p=p'qp''$, then it is clear that $|p'|,|p''|\ge d$. This proves that
$|w|\ge N+2d$.

Now let $p_0$ be a path of length $d$ that connects the identity with
some vertex $v$ from $\supp(\mu)$. We know from Lemma \ref{connflow} that
there is a loop $p$ at $v$ that consists of $N$ edges such that each
$P$-edge $e$ occurs in $p$ exactly $\mu(e)$ times. The path
$p_0pp_0^{-1}$ has length $N+2d$ and its label represents $g$. Therefore,
$|g|=N+2d$.
\endproof

Now we are ready to present an example of a strict dead end element
of depth $2$ in a free metabelian group. For simplicity, we consider
the group with $2$ generators. Let $\mu$ be a balanced flow given by
the ``thick" curve on the picture below. The curve or one of its
appropriate cyclic shifts read in the clockwise direction will be
denoted by $p$.

\begin{center}
\unitlength=1mm
\linethickness{0.8pt}
\begin{picture}(40.00,48.00)
\put(10.00,0.00){\line(0,1){10.00}}
\put(10.00,10.00){\line(-1,0){10.00}}
\put(0.00,10.00){\line(0,0){0.00}}
\put(0.00,10.00){\line(0,1){20.00}}
\put(0.00,30.00){\line(1,0){10.00}}
\put(10.00,30.00){\line(0,1){10.00}}
\put(10.00,40.00){\line(1,0){20.00}}
\put(30.00,40.00){\line(0,-1){10.00}}
\put(30.00,30.00){\line(1,0){10.00}}
\put(40.00,30.00){\line(0,-1){20.00}}
\put(40.00,10.00){\line(-1,0){10.00}}
\put(30.00,10.00){\line(0,-1){10.00}}
\put(30.00,0.00){\line(-1,0){20.00}}
\put(15.00,40.00){\vector(1,0){2.00}}
\put(10.00,36.00){\vector(0,0){0.00}}
\linethickness{0.4pt}
\put(10.00,10.00){\line(0,1){20.00}}
\put(10.00,30.00){\line(1,0){20.00}}
\put(30.00,30.00){\line(0,-1){20.00}}
\put(30.00,10.00){\line(-1,0){20.00}}
\put(0.00,20.00){\line(1,0){40.00}}
\put(20.00,0.00){\line(0,1){40.00}}
\put(20.00,20.00){\circle*{1.00}}
\put(18.00,17.00){\makebox(0,0)[cc]{$1$}}
\put(16.00,43.00){\makebox(0,0)[cc]{$a$}}
\put(7.00,35.00){\makebox(0,0)[cc]{$b$}}
\end{picture}

\nopagebreak
\end{center}

\noindent
It is easy to see that the flow $\mu$ determines an element
$g\in G=\fff_2/\fff_2''$, where
$g=b^2(ab^{-1})^2(b^{-1}a^{-1})^2(a^{-1}b)^2(ba)^2b^{-2}$. According
to Lemma \ref{lng}, the length of $g$ in $G$ equals $20$. In fact,
there are many ways to express $g$ by a word of minimal length.
One can take any path $q$ that connects the identity to a point
on the curve $p$, then go around the curve and return back to the
identity by the path $q^{-1}$. It is easy to see that the element $gx$
has length $19$ for any $x\in\{\,a^{\pm1},b^{\pm1}\,\}$. Indeed, one can
choose $q$ in such a way that the first edge $e$ of $q$ has label $x$.
If $q=eq'$, then the path $qp(q')^{-1}$ representing $gx$ has length $19$.
Therefore, the length of $gx$ equals $19$ (it cannot decrease by more
than $1$).

Now let $xy$ be any freely irreducible word of length $2$ over the
alphabet $\{\,a^{\pm1},b^{\pm1}\,\}$. Clearly, the path $q$ of
length $2$ labelled by $xy$ starting at the identity will end on the
curve $p$. Then $gxy$ is represented by the path $qp$ of length $18$,
which is also the length of $gxy$ in $G$.

By definition, $g\in G$ is a strict dead end element of depth $2$.
It is easy to see that $g$ is not a strict dead end of depth $3$
since $gaba^{-1}$ has length $19$.
\vspace{2ex}

A group word in the alphabet $\{x_1,\ldots,x_m\}$ is called
{\em positive\/} whenever it does not contain letters of the form
$x_i^{-1}$ ($1\le i\le m$). Here is one of the main results.

\begin{thm}
\label{oneend}
Let $R$ be a normal subgroup in the free group $\fff_m$ of rank
$m>1$. Assume that the following three conditions hold:

$1)$ the group $\fff_m/R$ has exactly one end,

$2)$ there are no nonempty positive words in $R$,

$3)$ any nonempty cyclically reduced word in $R$ has length
at least $\rho$.

Then the group $G=\fff_m/R'$ has strict dead end elements of depth
$[\rho/2]-1$.
\end{thm}

\proof
Let $k=[\rho/2]-1$. It follows from the third condition that the
ball $B_k$ of radius $k$ around the identity is a tree. Let $\ccc$
be the Cayley graph of $\fff_m/R$ and let $\ccc'$ be the subgraph
$\ccc-B_{k-1}$. We claim that $\ccc'$ is connected.

Assume the contrary. Then $\ccc'$ has more than one connected
component. Only one of them is infinite because of the first condition.
So there is a finite component $\Gamma$ in $\ccc'$. Each vertex in
$\ccc'$ has degree at least $2m-1>m$ since $m\ge2$. Let us take
any edge $e_1$ in $\Gamma$ with positive label and let $v_1$ be its
terminal vertex. There exists at least one edge $e_2$ with positive
label that comes out of $v_1$. (Otherwise $v_1$ has degree at most $m$.)
Denoting the terminal vertex of $e_2$ by $v_2$, we see that there is an
edge $e_3$ that comes out of $v_2$ and has a positive label, and so on.
Since $\Gamma$ is finite, we will have repetition of vertices on some
step and so there will be a positive nontrivial relation in $R$, which
contradicts the second condition.

Now we want to prove a slightly stronger fact. Suppose that $e$ is
an edge that belongs to $\ccc'$. We want to show that the graph
$\ccc'-e$ (we remove the edges $e^{\pm1}$) is still connected.
Assume the contrary. Suppose that the endpoints of $e$ belong to
different connected components of $\ccc'-e$. One of these components
must be finite. Let $v_0$ be the endpoint of $e$ that belongs to
this finite component $\Gamma$. All vertices in $\Gamma$ still have
degree at least $2m-1$ in this component except for, possibly, $v_0$,
which has degree at least $2m-2$. Let $e_1$ be an edge that comes out
of $v_0$. If this edge has a positive label, then we can repeat the
arguments from the previous paragraph. If the edge $e_1$ has a negative
label, then we repeat the same arguments replacing ``positive" by
``negative" everywhere and get a contradiction.

The fact that $\ccc'-e$ is always connected means that each edge
$e$ in $\ccc'$ belongs to a simple loop contained in $\ccc'$.
We want to show that $\ccc'$ contains a finite subgraph with certain
properties. First of all, we connect each pair of vertices at a
distance $k$ from the identity by a path in $\ccc'$. This gives a
finite subgraph in $\ccc'$. Now, for each edge of this subgraph,
we choose a simple loop in $\ccc'$ that contains this edge. We add
all these loops to the subgraph and get a new finite connected subgraph
$\Sigma$. By the construction, each edge of $\Sigma$ is contained in a
simple loop that is contained in $\Sigma$.

Let $e_1$, $e_2$, \dots, $e_t$ be all positively labelled edges of
$\Sigma$ and let $q_1$, $q_2$, \dots, $q_t$ be simple loops that
contain these edges, respectively (we allow some of the loops coincide).
For every $1\le i\le t$, let us consider the flow $2^i\chi_{q_i}$ and
let $\mu$ be the sum of these flows. Clearly, $\mu$ is balanced and
$\mu(e_i)\ne 0$ for all $1\le i\le t$. Indeed, for any $e_i$ let $j$ be
the smallest number such that $e_i^{\pm1}$ occurs in $q_j$. Clearly,
$j\le i$ and $\mu(e_i)$ equals $2^j$ modulo $2^{j+1}$.

The identity element is at a distance $k$ from $\Sigma$. If $N$ is
the weight of $\mu$, then the element $g$ in $\fff_m/R'$ represented
by $\mu$ has length $N+2k$ in this group according to Lemma \ref{lng}.
We claim that $g$ is a strict dead end element of depth at least $k$
in $\fff_m/R'$.

Indeed, let $a_1\cdots a_k$ be any freely irreducible word in the
generators of $G$. We want to prove that
$|g|>|ga_1|>|ga_1a_2|>\cdots>|ga_1a_2\cdots a_k|$. We show that
$|ga_1\cdots a_i|=|g|-i$ for all $0\le i\le k$. Obviously,
$|ga_1\cdots a_{i-1}a_i|\ge|g|-i$. On the other hand, let $q$ be the
path of length $k$ that connects $1$ to a vertex $v$ in $\Sigma$
such that $q=q'q''$ and the label of $q'$ is $a_1a_2\cdots a_i$. By
Lemma \ref{connflow}, there exists a loop $p$ at $v$ of length $N$ that
induces $\mu$. Then the label of the path $qp(q'')^{-1}$ represents
$ga_1\cdots a_i$. The length of this path is $|g|-|q'|=|g|-i$. This
implies what we wanted to prove.
\endproof

Let us estimate the length $\rho_d$ of a shortest nontrivial relation
in $\Sol(m,d)$, where $m>1$. It is not hard to see that $\rho_1=4$,
$\rho_2=14$. In general, we have the following rough estimate.

\begin{lm}
\label{shortrel}
Let $\rho_d$ be the length of the shortest nontrivial relation
in the free soluble group $\Sol(m,d)$, where $m>1$, $d\ge1$. Then
$\rho_{d+1}\ge3\rho_d$ for all $d\ge1$.
\end{lm}

\proof
Let $R=\fff_m^{(d)}$. Suppose that $w$ is a shortest nontrivial
word that represents $1$ in $\fff_m/R'$. By Lemma \ref{wp}, $w$
induces a zero flow on the Cayley graph $\ccc$ of $\fff_m/R$. Let
$p=p(w)$ be the path in $\ccc$ labelled by $w$ starting at $1$.
Clearly, $p$ cannot have (cyclic) subpaths of the form $ee^{-1}$,
where $e$ is an edge. (Otherwise $w$ is not a shortest nontrivial
relation.)

Let $\Gamma$ be a subgraph in $\ccc$ formed by all the edges that
occur in $p$. It is obvious that $\Gamma$ is connected; it is also
clear that $\Gamma$ has no vertices of degree $1$.

One can observe that $\Gamma$ is not a circle, that is, it cannot
have rank $1$ (a rank of a connected graph is the rank of its fundamental
group). Hence the rank of $\Gamma$ is at least $2$. It is well known that
each graph of rank greater than $1$ has a (geometric) subgraph of one of
the following forms (cf. \cite{GT}):

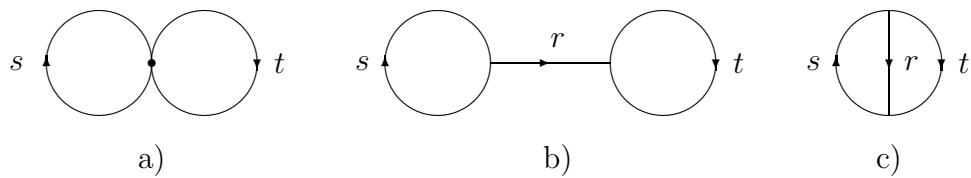
\begin{center}
\unitlength=1.00mm
\linethickness{0.4pt}
\begin{picture}(127.00,21.00)
\put(12.00,14.00){\circle{14.00}}
\put(26.00,14.00){\circle{14.00}}
\put(19.00,14.00){\circle*{1.00}}
\put(57.00,14.00){\circle{14.00}}
\put(87.00,14.00){\circle{14.00}}
\put(64.00,14.00){\line(1,0){16.00}}
\put(117.00,14.00){\circle{14.00}}
\put(117.00,7.00){\line(0,1){14.00}}
\put(5.00,14.00){\vector(0,1){1.00}}
\put(33.00,14.00){\vector(0,-1){1.00}}
\put(50.00,14.00){\vector(0,1){1.00}}
\put(71.00,14.00){\vector(1,0){1.00}}
\put(94.00,14.00){\vector(0,-1){1.00}}
\put(110.00,14.00){\vector(0,1){1.00}}
\put(117.00,14.00){\vector(0,-1){1.00}}
\put(124.00,14.00){\vector(0,-1){1.00}}
\put(1.00,14.00){\makebox(0,0)[cc]{$s$}}
\put(36.00,14.00){\makebox(0,0)[cc]{$t$}}
\put(47.00,14.00){\makebox(0,0)[cc]{$s$}}
\put(97.00,14.00){\makebox(0,0)[cc]{$t$}}
\put(73.00,17.00){\makebox(0,0)[cc]{$r$}}
\put(107.00,14.00){\makebox(0,0)[cc]{$s$}}
\put(127.00,14.00){\makebox(0,0)[cc]{$t$}}
\put(120.00,14.00){\makebox(0,0)[cc]{$r$}}
\put(19.00,1.00){\makebox(0,0)[cc]{a)}}
\put(73.00,1.00){\makebox(0,0)[cc]{b)}}
\put(117.00,1.00){\makebox(0,0)[cc]{c)}}
\end{picture}

\nopagebreak
\end{center}

We see that a) can be considered as a partial case of b) with $r$
as the empty path. Since $p$ goes at least once through each edge
of $\Gamma$ in both directions, one can estimate the length of $p$
as $|p|\ge2(|s|+|t|+|r|)$ in b) and c). Since every nontrivial
relation in $\fff_m/R$ has length at least $\rho_d$, we can conclude
that $|s|,|t|\ge r_d$ in the case b) and so $|p|\ge4\rho_d$. In the
case c), we have loops $sr$, $r^{-1}t$, and $st$, each of length
at least $r_d$. Hence $|p|\ge2(|s|+|t|+|q|)=(|s|+|q|)+(|q|+|t|)+
(|s|+|t|)\ge3\rho_d$. Therefore, we always have $\rho_{d+1}\ge3\rho_d$.
\endproof

Lemma \ref{shortrel} immediately implies that $\rho_d\ge4\cdot3^{d-1}$
for the group $\Sol(m,d)$, where $m>1$, $d\ge1$.

Now we can apply Theorem \ref{oneend} to extract a corollary for
the case of free soluble groups.

\begin{thm}
\label{freesol}
If $m>1$, $d\ge2$, then the free soluble group $\Sol(m,d)$ has
strict dead end elements of depth at least $2\cdot3^{d-2}-1$.
\end{thm}

\proof
Free soluble groups have only one end in the sense of Stallings.
Indeed, they are not virtually cyclic and they cannot be presented
as a nontrivial amalgamated product or an HNN-extension. Let
$R=\fff_m^{(d-1)}$. Every nontrivial word in the commutator subgroup
of $\fff_m$ has zero exponent sum on each generator. Hence the only
positive word in $R\subseteq\fff_m'$ is the empty word. It follows
from Lemma \ref{shortrel} that the length of the shortest nontrivial
relation in $R$ is $\rho_{d-1}\ge4\cdot3^{d-2}$. Now we apply
Theorem \ref{oneend} and conclude that the group $\Sol(m,d)=\fff_m/R'$
has strict dead end elements of depth at least
$[\,\rho_{d-1}/2\,]-1\ge2\cdot3^{d-2}-1$.
\endproof

Notice that the groups $\Sol(m,d)$ are not finitely presented whenever
$m>1$, $d\ge2$ by a result of \cite{Shm}. However, they are residually
finite \cite{Gru}. Therefore, for any $k\ge1$ there exist finitely
presented and even finite groups with strict dead end elements of
depth $k$.
\vspace{2ex}

{\bf Remark.}\ The solution to the word problem in free soluble
groups gives rise to a question about growth functions for these
groups. Recall that if $\ccc=\ccc(G,A)$, then the {\em growth function\/}
of $G$ with respect to $A$ is the sequence $b_n=\Card(B_n)$, where $B_n$
denotes the ball of radius $n$ in $\ccc$. The limit
$\beta=\lim_{n\to\infty}\sqrt[n]{b_n}$ always exists; it is called the
{\em growth rate\/} of $G$ (with respect to $A$).

Even in the case of $G=\Sol(2,2)$, a $2$-generated free metabelian
group, neither the growth function nor the growth rate is known.
A lower bound for the growth function of $G$ is the number of
self-avoiding walks (SAW) of length $n$ in the square lattice starting
at the identity (different self-avoiding paths represent different
elements of $G$ by Lemma \ref{wp}). The exact value of the growth
rate of SAW is not known although more than 10 of its decimal digits
were found \cite{GC}. So the growth rate of $G$ satisfies inequality
$\beta\ge2.63815853034$.

Notice that for some metabelian $2$-generated groups both the growth
rate and the growth function are known. Say, Parry \cite{Par} found
them for the restricted wreath product $\zz\wr\zz$ and for some other
groups. For $\zz\wr\zz$, the growth function is rational and the growth
rate equals $\sqrt{2}+1$. In \cite{CEG} it is shown that the growth
functions of Baumslag -- Solitar groups $BS(1,n)$ are rational; their
growth rates do not exceed the above value $\sqrt{2}+1$.

In a recent paper \cite{AGG}, the authors prove that the growth
rate of the free soluble group $\Sol(m,d)$ approaches $2m-1$ as
$d\to\infty$. This shows that the growth rate of amenable groups
with $m$ generators can be arbitrarily close to the maximum value
$2m-1$.

\end{document}